\documentclass[12pt]{amsart}
\usepackage{amscd,amsmath,amssymb,amsfonts}
\theoremstyle{plain}
\newtheorem{thm}{Theorem}
\newtheorem{lem}[thm]{Lemma}
\newtheorem{cor}[thm]{Corollary}
\newtheorem{prop}[thm]{Proposition}
\theoremstyle{definition}
\newtheorem{remark}[thm]{Remark}
\newtheorem{defn}[thm]{Definition}
\newtheorem{ex}[thm]{Example}
\newtheorem{conj}[thm]{Conjecture}
\numberwithin{thm}{section}
\numberwithin{equation}{section}
\newcommand{\eq}[2]{\begin{equation}\label{#1}#2 \end{equation}}
\newcommand{\ml}[2]{\begin{multline}\label{#1}#2 \end{multline}}
\newcommand{\ga}[2]{\begin{gather}\label{#1}#2 \end{gather}}

\newcommand{\inj}{\hookrightarrow}

\newcommand{\Spec}{{\rm Spec \,}}

\newcommand{\tr}{{\rm Tr}}
\newcommand{\eps}{{\epsilon}}
\newcommand{\sA}{{\mathcal A}}

\newcommand{\sC}{{\mathcal C}}

\newcommand{\sG}{{\mathcal G}}
\newcommand{\sH}{{\mathcal H}}

\newcommand{\sL}{{\mathcal L}}
\newcommand{\sM}{{\mathcal M}}
\newcommand{\sN}{{\mathcal N}}
\newcommand{\sO}{{\mathcal O}}
\newcommand{\sP}{{\mathcal P}}

\newcommand{\sZ}{{\mathcal Z}}
\newcommand{\A}{{\mathbb A}}

\newcommand{\C}{{\mathbb C}}

\newcommand{\F}{{\mathbb F}}
\newcommand{\G}{{\mathbb G}}
\renewcommand{\H}{{\mathbb H}}

\renewcommand{\P}{{\mathbb P}}
\newcommand{\Q}{{\mathbb Q}}

\newcommand{\V}{{\mathbb V}}

\newcommand{\Z}{{\mathbb Z}}

\begin{document}

\title{The Additive Dilogarithm}
\author{Spencer Bloch}
\address{Dept. of Mathematics,
University of Chicago,
Chicago, IL 60637,
USA}
\email{bloch@math.uchicago.edu}

\author{H\'el\`ene Esnault}
\address{Mathematik,
Universit\"at Essen,
Essen,
Germany}
\email{esnault@uni-essen.de}
\date{Oct. 8, 2002}
\dedicatory{To Kazuya Kato, with fondness and profound respect, on the
occasion of his fiftieth birthday}
\thanks{ This work has been partly supported by the NSF grant
DMS-9423007-A1, the DFG Schwerpunkt ''Komplexe Geometrie'',
the Humboldt Foundation, the University Paris VII, the IHES and
the Newton Institute.}

\maketitle

\section{Introduction}

In \cite{Lau} D\'ef. (5.1.1), Laumon introduces the category of
generalized $1$-motives over a field $k$ of characteristic $0$. Objects
in this category are arrows
$f:\sG \to G$ where
$\sG$ and $G$ are commutative algebraic groups, with $\sG$ assumed
formal, torsion free, and $G$ assumed connected. These, of course,
generalize the more restricted category of $1$-motives introduced by
Deligne \cite{D} as a model for the category of mixed Hodge structures of
types $\{(0,0), (0,-1), (-1,0), (-1,-1)\}$. Of particular interest for
us are motives of the form $\Z \to \V$ which arise in the study of
algebraic cycles relative to a ``modulus''. Here $\V\cong \G_a^n$ is a
vector group. The simplest example is
\eq{1.1}{\text{Pic}(\A^1,2\{0\}) \cong \G_a.
}
which may be viewed as a degenerate version of the identification
$\text{Pic}(\A^1,\{0,\infty\}) \cong \G_m$ obtained by associating to a
unit the corresponding Kummer extension of $\Z$ by $\Z(1)$. (For more
details, cf.
\cite{BE},\cite{BS},\cite{E},\cite{ESV},\cite{L}.) We expect such
generalized motives to play an important role in the (as yet undefined)
contravariant theory of motivic sheaves and motivic cohomology for
(possibly singular) varieties.

The polylog mixed motives of Beilinson and Deligne are generalizations
to higher weight of Kummer extensions, so it seems natural to look for
degenerate, or $\G_a$ versions of these. The purpose of this article is to
begin to study an additive version of the dilogarithm motive. We
assume throughout that $k$ is a field which for the most part will be
taken to be of characteristic
$0$. Though our results are limited to the dilogarithm,
the basic result from cyclic homology
\eq{1.2}{Gr^n_\gamma\ker\Big(K_{2n-1}(k[t]/(t^2)) \to K_{2n-1}(k)\Big)
\cong k
}
suggests that higher polylogarithms exist as well.

In the first part of the article we introduce an additive
``Bloch group'' $TB_2(k)$ for an algebraically closed field $k$ of
characteristic $\neq 2$. In lieu of the $4$-term sequence in
motivic cohomology associated to the usual Bloch group \ml{1.3}{0
\to H^1_M(\Spec(k),\Q(2)) \to B_2(k) \to k^\times \otimes
k^\times\otimes \Q \\
\to H^2_M(\Spec(k),\Q(2)) \to 0
}
(with $H^1_M(\Spec(k),\Q(2)) \cong K_3(k)_{ind} \otimes \Q$ and
$H^2_M(\Spec(k),\Q(2)) \cong K_2(k)\otimes \Q$), we find an additive
$4$-term sequence
\ml{1.4}{0 \to TH^1_M(\Spec(k),\Q(2)) \to TB_2(k) \to k \otimes
k^\times \\
\stackrel{d\log}{\to} TH^2_M(\Spec(k),\Q(2)) \to 0
}
where
\ga{1.5}{TH^1_M(\Spec(k),\Q(2)) := K_2(\A^1_t,(t^2)) \cong
(t^3)/(t^4)\cong k;\\
TH^2_M(\Spec(k),\Q(2)) := K_1(\A^1,(t^2)) \cong \Omega^1_k =
\text{absolute K\"ahler
$1$-forms};\notag \\
d\log(a\otimes b) = a\frac{db}{b}. \notag }

Our construction should be compared and contrasted with the results of
\cite{C}. Cathelineau's group $\beta_2(k)$ is simply the kernel
\eq{1.6}{0\to \beta_2(k) \to k\otimes k^\times \to \Omega^1_k \to 0,
}
so there is an exact sequence \minCDarrowwidth.1cm
\eq{1.7}{\begin{CD}0 @>>> TH^1_M(\Spec(k),\Q(2)) @>>> TB_2(k) @>>>
\beta_2(k)
\to 0 \\
@. @VV \cong V \\
@. k \end{CD} }
For $a\in k$ we define $\langle a\rangle \in TB_2(k)$ lifting similar
elements defined by Cathelineau and satisfying his $4$-term
infinitesimal version
\eq{1.8}{\langle a\rangle - \langle b\rangle + a\langle b/a\rangle +
(1-a)\langle (1-b)/(1-a)\rangle = 0;\quad a\neq 0,1.
}
of the classical $5$-term dilogarithm relation.
Here, the notation $x\langle y\rangle$ refers to an action of $k^\times$
on $TB_2(k)$. Unlike $\beta_2(k)$, this action does {\it not} extend to a
$k$-vector space structure on $TB_2(k)$. Thus \eqref{1.7} is an exact
sequence of $k^\times$-modules, where the kernel and cokernel have
$k$-vector space structures but the middle group does not.

Finally in this section we show the assignment $\langle a\rangle \mapsto
a(1-a)$ defines a regulator map $\rho : TB_2(k) \to k$ and the composition
\eq{1.9}{TH^1_M(\Spec(k),\Q(2)) \inj TB_2(k) \stackrel{\rho}{\to} k
}
is an isomorphism.

It seems plausible that $TB_2(k)$ can be interpreted as a Euclidean
sissors-congruence group, with $\partial: TB_2(k) \to k\otimes k^\times$
the Dehn invariant and $\rho:TB_2(k) \to k$ the volume. Note the scaling
for the $k^\times$-action is appropriate, with $\partial(x\langle y\rangle)
= x\partial(\langle y\rangle)$ and $\rho(x\langle y\rangle) =
x^3\rho(\langle y\rangle)$. For a careful discussion of Euclidean
sissors-congruence
and its relation with the dual numbers, the reader is referred to \cite{Gr}
and the references cited there.

In \S 4 we introduce an extended polylogarithm Lie algebra. The dual co-Lie
algebra has generators
$\{x\}_n$ and $\langle x\rangle_n$ for $x\in k-\{0,1\}$. The dual of the
bracket satisfies $\partial \{x\}_n = \{x\}_{n-1}\cdot\{1-x\}_1$ and
$\partial \langle x\rangle_n = \langle x\rangle_{n-1}\cdot \{1-x\}_1 +
\langle 1-x\rangle_1\cdot \{x\}_{n-1}$ with $\langle x\rangle_1 = x \in
k$. For example, $\partial \langle x\rangle_2 = x\otimes x + (1-x) \otimes
(1-x) \in k\otimes k^\times$ is the Cathelineau relation \cite{C}. It seems
likely that there exists a representation of this Lie algebra, extending the
polylog
representation of the sub Lie algebra generated by the $\{x\}_n$, and
related to variations of Hodge structure over the dual numbers lifting the
polylog Hodge structure.

\S 5 was inspired by Deligne's interpretation of symbols \cite{Dsymb} in
terms of line
bundles with connections. We indicate how this viewpoint is related to the
additive dilogarithm. In characteristic $0$, one finds affine
bundles with connection, and the regulator map on $K_2$ linearizes to the
evident map $H^0(X,\Omega^1) \to \H^1(X, \sO \to \Omega^1)$. In
characteristic $p$, Artin-Schreier yields an exotic flat realization of the
additive dilogarithm motive. For simplicity we limit ourselves to
calculations mod $p$. The result is a flat covering $T$ of $\A^1-\{0,1\}$
which is a torsor under a flat Heisenberg groupscheme $\sH_{AS}$. This
groupscheme has a natural representation on the abelian groupscheme $\V :=
\Z/p\Z \oplus \mu_p\oplus \mu_p$. The contraction
\eq{1.10a}{ T \stackrel{\sH_{AS}}{\times} \V
}
should, we think, be considered as analogous to the mod $\ell$ \'etale
sheaf on $\A^1-\{0,1\}$ with fibre $\Z/\ell\Z \oplus \mu_\ell \oplus
\mu_\ell^{\otimes 2}$ associated to the $\ell$-adic dilogarithm.

The polylogarithms can be interpreted in terms of algebraic
cycles on products of copies of $\P^1-\{1\}$ (\cite{B}, (3.3)),
so it seems natural to consider algebraic cycles on
\eq{1.10}{(\A^1,2\{0\})\times(\P^1-\{1\},\{0,\infty\})^{n}. } In
the final section of this paper, we calculate the Chow groups of
$0$-cycles on these spaces. Our result:
\eq{1.11}{CH_0\Big((\A^1,2\{0\})\times(\P^1-\{1\},\{0,\infty\})^{n}\Big)
\cong \Omega^n_k, \ \ n\ge 0, } is a ``degeneration'' of the
result of Totaro \cite{T} and Nesterenko-Suslin \cite{NS}
\eq{1.12}{CH_0\Big((\P^1-\{1\},\{0,\infty\})^{n}\Big) \cong
K_n^M(k) = \text{$n$-th Milnor $K$-group,} } and a cubical
version of the simplicial result $SH^n(k,n)\cong \Omega^{n-1}_k$
(see \cite{BE}).

We thank J\"org Wildeshaus and Jean-Guillaume Grebet for helpful remarks.

\section{Additive Bloch groups}

Let $k$ be a field with $1/2 \in k$. In this section, we
mimic the construction in \cite{Bbook} \S 5, replacing
the semi-local ring of functions on $\P^1$, regular at $0$ and
$\infty$ by the local ring of functions on $\A^1$, regular at
$0$, and the relative condition on $K$-theory at $0$ and $\infty$
by the one at $2\cdot\{0\}$. In particular, as we fix only $0$ and
 $\infty$ in this theory, we have a $k^\times$-action on the parameter
 $t$ on $\A^1$ so our groups will be $k^\times$-modules.

 Thus let $R$ be the local ring at $0$ on $\A^1_k$.
One has an exact sequence of relative $K$-groups
\eq{1}{K_2(\A^1_k) \to K_2(k[t]/(t^2)) \to K_1(\A^1,(t^2)) \to
K_1(\A^1) \to K_1(k[t]/(t^2)) .}
Using Van der Kallen's
calculation of $K_2(k[t]/(t^2))$ \cite{VK} and the homotopy
property $K_*(k) \cong K_*(\A^1_k)$, we conclude
\eq{2}{K_1(\A^1_k,(t^2)) \cong \Omega^1_k. } Now we localize on
$\A^1$ away from $0$. Assuming for simplicity that $k$ is
algebraically closed, we get \eq{3}{\coprod_{k-\{0\}} K_2(k) \to
K_2(\A^1,(t^2)) \to K_2(R,(t^2)) \to \coprod_{k-\{0\}} k^\times
\to \Omega^1_k \to 0 }

To $a\in (t^2)$ and $b\in R$ we associate the pointy-bracket
symbol \cite{Lo} $\langle a,b\rangle\in K_2(R,(t^2))$ which
corresponds to the Milnor symbol $\{1-ab,b\}$ if $\ b\neq 0$.
These symbols generate $K_2(R, (t^2))$. If the divisors of
$a$ and $b$ are disjoint, we get \eq{4}{\text{tame}\langle
a,b\rangle = a|_{\text{poles of }b} + b|_{ab=1} +
b^{-1}|_{\text{poles of }a} }
We continue to assume $k$
algebraically closed. Let $\sC \subset K_2(R,(t^2))$ be the
subgroup generated by pointy-bracket symbols with $b\in k$. For
$a\in (t^2)$ write
\eq{2.5}{a(t)=\frac{a_0t^n+\ldots+a_{n-2}t^2}{t^m+b_1t^{m-1}+\ldots
+b_{m-1}t+b_m};\quad b_m \neq 0. } We assume numerator and
denominator have no common factors.
If $\alpha_i$ are the solutions to the equation $a(t) =
\kappa\in k^\times\cup \infty$, then $\sum \alpha_i^{-1} =
-b_{m-1}/b_m$. In particular, this is independent of
$\kappa$. It follows that one has an isomorphism
\eq{6}{\coprod_{k-\{0\}}k^\times\Big/\text{tame}(\sC) \cong
k\otimes_\Z k^\times;\quad u|_v \mapsto v^{-1}\otimes u. }

Define
\ga{7}{TB_2(k) := K_2(R,(t^2))/\sC \\
TH^1_M(k,2) := \text{image}
\Big(K_2(\A^1,(t^2)) \to TB_2(k)\Big) \notag \\
TH^2_M(k,2) := \Omega^1_k = K_1(R,(t^2))\notag
}
A basic result of Goodwillie \cite{Good} yields $K_2(\A^1,(t^2))\cong k$,
so $TH^1_M(k,2)$ is a quotient of $k$. We will see (remark \ref{rem2.5})
that in fact $TH^1_M(k,2)  \cong k$.
The above discussion yields
\begin{prop}\label{prop2.1} Let $k$ be an algebraically closed field of
  characteristic
  $\neq 2$. With notations as above, we have an exact sequence
\eq{8}{0 \to TH^1_M(k,2) \to TB_2(k) \stackrel{\partial}{\to} k\otimes
  k^\times \stackrel{\pi}{\to}
\Omega^1_k \to 0.
}
Here $\pi(a\otimes b) = a\frac{db}{b}$ and $\partial$ is defined via the
  tame symbol. \end{prop}
\begin{remark}\label{rem2.2} There is an evident action of the group
  $k^\times$ on $\A^1$ (multiplying the parameter) and
  hence on the sequence \eqref{8}. This action extends to a $k$-vector
  space structure on all the terms except $TB_2(k)$.
\end{remark}

Let $\mathfrak m = tR\subset R$. One has the following purely
algebraic description of $K_2(R,\mathfrak m^2)$ (\cite{S}, formula
(1.4), and the references cited there). \newline \noindent
generators: \eq{2.9}{\langle a,b \rangle; \qquad (a,b) \in
(R\times \mathfrak m^2) \cup (\mathfrak m^2 \times R) } \noindent
Relations: \ga{10}{ \langle a,b \rangle = -\langle b,a
\rangle;\quad a \in \mathfrak m^2
\\
\langle a,b \rangle + \langle a,c \rangle = \langle a,b+c-abc
\rangle;\quad a\in \mathfrak m^2 \text{ or } b,c \in \mathfrak m^2 \\
\langle a,bc \rangle = \langle ab,c \rangle + \langle ac,b \rangle; \quad
a \in \mathfrak m^2
 }
\begin{prop}\label{prop2} There is a well-defined and nonzero map
\eq{11}{ \rho:K_2(R,\mathfrak m^2)
\to \mathfrak m^3/\mathfrak m^4
}
defined by
\eq{12}{\rho\langle a,b \rangle := \begin{cases} -adb & a\in \mathfrak m^2 \\
bda & b\in \mathfrak m^2.
\end{cases}
}
\end{prop}
\begin{proof} Note first if $a,b\in \mathfrak m^2$ then $adb\equiv bda \equiv
0 \mod \mathfrak m^4$ so the definition \eqref{12} is consistent. For $a\in
\mathfrak m^2$
\eq{13}{\langle a,b \rangle + \langle b,a \rangle \mapsto -adb+adb=0,
}
so \eqref{10} holds. For $a \in \mathfrak m^2$
\eq{14}{\langle a,b \rangle + \langle a,c \rangle \mapsto -ad(b+c) \equiv
-ad(b+c-abc) \mod \mathfrak m^4
}
for $b,c \in \mathfrak m^2$
\eq{15}{\langle a,b \rangle + \langle a,c \rangle \mapsto (b+c)da \equiv
(b+c-abc)da \mod \mathfrak m^4
}
For $a\in \mathfrak m^2$,
\eq{16}{\langle a,bc \rangle \mapsto -ad(bc) = -abdc-acdb = \rho(\langle
ab,c \rangle+\langle ac,b \rangle)
}
\end{proof}

\begin{remark} Note that \eq{117}{-adb \equiv \log(1-ab)db/b \in
\mathfrak{m}^2\Omega^1_R/d\log (1+\mathfrak{m}^4)\cong
\mathfrak{m}^3/\mathfrak{m}^4. } The group $\mathfrak{m}^2\Omega^1_R/d\log
(1+\mathfrak{m}^4)$ is the group of isomorphism classes of rank 1
line bundles, trivialized at the order 4 at $\{0\}$, with a
connection vanishing at the order 2 at $\{0\}$. Thus the
regulator map $\rho$ assigns such a connection to a pointy symbol.
Over the field of complex numbers $\C$, one can think of it in
terms of ``Deligne cohomology'' $\H^2(\A^1, j_!\Z(2)\to t^4 \sO
\to t^2\omega)$, and one can, as in \cite{ELoday}, write down
explicitely an analytic Cech cocycle for this regulator as a
Loday symbol.
\end{remark}

 One has $\rho\langle t^2,x \rangle = -t^2dt
\neq 0$, thus $\rho $ is not trivial. Note also, the
appearance of $\mathfrak m^3$ is consistent with  A. Goncharov's
idea \cite{G} that the regulator in this context should
correspond to the volume of a simplex in hyperbolic $3$ space in
the sissors-congruence interpretation \cite{Gr}. In particular, it
should scale as the third power of the coordinate.

 Proposition
\ref{prop2} yields
\begin{cor}\label{cor2.4} Let $\mathfrak m\subset R$ be the maximal
ideal. One has a well-defined map
\eq{19}{\rho: TB_2(k) \to \mathfrak m^3/\mathfrak m^4
}
given on pointy-bracket symbols by
\eq{20}{\rho\langle a,b\rangle = -a\cdot db;\quad a\in \mathfrak m^2,\
b\in R.
}
For $x\in TB_2(k)$ and $c\in k^\times$, write $c\star x$ for the image of
$x$ under the mapping $t\mapsto c\cdot t$ on polynomials. Then $\rho(c\star
x) = c^3\cdot\rho(x)$.
\end{cor}
\begin{proof} The first assertion follows because if $b\in k$, then $db=0$.
The second assertion is clear.
\end{proof}

\begin{remark}\label{rem2.5} The map $\rho$ is non-trivial on
  $TH^1_M(k,2)$ because
$\rho\langle t^2,t \rangle = -t^2dt \neq 0$. Since this group is a
  $k^\times$-module (remark \ref{rem2.2}) and is a quotient of $k$ by the
  result of Goodwillie cited above, it follows that
\eq{21}{TH^1_M(k,2) \cong (t^3)/(t^4) \cong k
}
\end{remark}

\section{Cathelineau elements and the entropy functional equation}

We continue to assume $k$ is an algebraically closed field of
characteristic $\neq 2$. Define for $a\in k-\{0,1\}$
\ga{3.1}{ \langle a\rangle := \langle t^2, \frac{a(1-a)}{t-1}\rangle \in
TB_2(k) \\
\eps(a) := a\otimes a + (1-a) \otimes (1-a) \in k^\times \otimes k \notag
}
\begin{lem}Writing $\partial$ for the tame symbol as in proposition
  \ref{prop2.1}, we have $\partial(\langle a \rangle) =2\eps(a)$.
\end{lem}
\begin{proof}
\ml{3.2a}{\partial(\langle a \rangle) =
\text{tame}\Big\{\frac{\frac{1-t}{a(1-a)}+t^2}{ \frac{1-t}{a(1-a)}},
\frac{a(1-a)}{t-1} \Big\} = \\
 \frac{a(1-a)}{t-1}\Big |_{t=\frac{1}{a}} +
\frac{a(1-a)}{t-1}\Big |_{t=\frac{1}{1-a}} \mapsto a^2\otimes a +(1-a)^2
\otimes
(1-a)  = 2\eps(a) \in k^\times \otimes k.
}
\end{proof}

\begin{lem}\label{lem3.2} We have $\rho(\langle a \rangle) = a(1-a)t^2dt
\in (t^3)/(t^4)$.
\end{lem}
\begin{proof} Straightforward from corollary \ref{cor2.4}.
\end{proof}

\begin{lem}\label{lem3.3}
Let notations be as in corollary \ref{cor2.4},
Assume $k$ is algebraically closed, and
${\rm char}(k)\neq 2, 3$. Then every element in $TB_2(k)$ can be written as a
sum $\sum c_i\star \langle a_i\rangle$. In other words, $TB_2(k)$
is generated as a $k^\times$-module by the $\langle a\rangle$.
\end{lem}
\begin{proof}Define
\eq{3.2}{\mathfrak b := \text{Image}(\partial:TB_2(k) \to
k^\times \otimes k) = \ker(k^\times \otimes k \to \Omega^1_k). }
The $k$-vectorspace structure  $c\cdot (a\otimes b)$ on
$k^\times \otimes k$ is defined by $a\otimes c b$.
By \eqref{6} and \eqref{4}, the map $TB_2(k)\to k^\times \otimes k$
is $k^\times$-equivariant.

Let $A\subset TB_2(k)$ be the subgroup generated by the $c\star
\langle a\rangle$. $\mathfrak b$ is a $k$-vector space which is
generated \cite{C} by the $\eps(a)$ so the composition $A\subset
TB_2(k) \to \mathfrak b$ is surjective.  For $c_1, c_2 \in
k^\times$ with $c_1+c_2\neq 0$ we have $(c_1+c_2)\star \langle
a\rangle -c_1\star \langle a\rangle - c_2\star \langle a\rangle
\mapsto 0 \in \mathfrak b$, so this element lies in $A\cap
H^1_M(k,2)$. It is not trivial because
\ml{3.3}{\rho((c_1+c_2)\star \langle a\rangle -c_1\star \langle
a\rangle -
c_2\star \langle a\rangle) = \\
\Big((c_1+c_2)^3 -c_1^3-c_2^3\Big)a(1-a)t^2dt =\\
3\Big(c_1c_2(c_1+c_2)\Big)a(1-a)t^2dt. } Since the equation
$\lambda=3\Big(c_1c_2(c_1+c_2)\Big)a(1-a)$ can be solved in $k$,
one has $A\supset H^1_M(k,2)$. This finishes the proof.
\end{proof}

\begin{thm}\label{thm3.4} Under the assumptions of
lemma \ref{lem3.3},  the group $TB_2(k)$ is generated as a
$k^\times$-module by the $\langle a\rangle$. These satisfy relations
\eq{3.4}{\langle a\rangle - \langle b\rangle + a\star\langle
b/a\rangle +(1-a)\star\langle (1-b)/(1-a)\rangle = 0.  
}
\end{thm}
\begin{proof}The generation statement is lemma \ref{lem3.3}. Because we
factor out by symbols with one entry constant, we get
\eq{3.5}{x\star\langle a\rangle = \langle x^2t^2,\frac{a(1-a)}{xt-1}
\rangle =
\langle t^2,\frac{x^2a(1-a)}{xt-1}\rangle.
}
The identity to be established then reads
\ml{3.6}{0= \langle t^2,\frac{a(1-a)}{t-1}\rangle - \langle
t^2,\frac{b(1-b)}{t-1}\rangle +
\langle t^2,\frac{b(a-b)}{at-1}\rangle +
\langle t^2,\frac{(1-b)(b-a)}{(1-a)t-1}\rangle.
}
The pointy bracket identity $\langle a,b\rangle + \langle a,c\rangle =
\langle a,b+c-abc\rangle$ means we can compute the above sum using
``faux'' symbols
\ml{3.8a}{\{t^2,1-\frac{a(1-a)t^2}{t-1}\}\{t^2,1-\frac{b(1-b)t^2}{t-1}\}^{-1}
\{t^2,1-\frac{b(a-b)t^2}{at-1}\}\times \\
\{t^2,1-\frac{(1-b)(b-a)t^2}{(1-a)t-1}\} = \{t^2,X\}
}
with
\ml{3.7}{X = \\
\frac{(1-t+a(1-a)t^2)(1-at+b(a-b)t^2)(1-(1-a)t+(1-b)(b-a)t^2)}{
(1-t+b(1-b)t^2)(1-at)(1-(1-a)t)} \\
=
\frac{(1-at)(1-(1-a)t)(1-bt)(1-(a-b)t)(1-(1-b)t)(1-(b-a)t)}{
(1-bt)(1-(1-b)t)(1-at)(1-(1-a)t)} \\
=(1-(a-b)t)(1-(b-a)t) = 1-(a-b)^2t^2
}
Reverting to pointy brackets, the Cathelineau relation equals
\eq{3.8}{\{t^2,X\} = \{1-(a-b)^2t^2,(a-b)^2\} = \langle t^2,(a-b)^2\rangle
= 0
}
since we have killed symbols with one entry constant.
\end{proof}
\begin{remark}One can get a presentation for $TB_2(k)$ if one imposes
\eqref{3.4} and in addition relations of the form
\ga{}{\Big((x+y+z+w)-(x+y+z)-(x+y+w)-\ldots -(x)-(y)-(z)-(w)\Big)\star \langle
a \rangle = 0 \\
(-1)\star \langle a\rangle = -\langle a\rangle.
}
Here $x\in k^\times$ corresponds to $(x) \in \Z[k^\times]$, and the first
relation is imposed whenever it makes sense, i.e. whenever all the partial
sums are non-zero. The proof uses uniqueness of solutions for the entropy
equation \cite{Da}. Details are left for the reader.
\end{remark}
\begin{remark}It
is remarkable that a functional equation equivalent to \eqref{3.4},
\eq{3.9}{
\langle a\rangle  + (1-a)\star\langle \frac{b}{1-a}\rangle =
\langle b\rangle + (1-b)\star\langle \frac{a}{1-b}\rangle }
occurs in information theory, where it is known to have a unique
continuous functional solution (up to scale) given by
$y\star\langle x\rangle \mapsto -yx\log(x)-y(1-x)\log(1-x)$. If on
the other hand, we interpret the torus action $y\star$ as
multiplication by $y^p, p\neq 1$, then the unique solution is
$\langle x\rangle \mapsto x^p+(1-x)^p -1$ \cite{Da}. Note the
regulator map $\rho(y\star\langle x\rangle) = y^3x(1-x)$, so
$\rho$ is a solution for $p=3$. Indeed, $x(1-x) =
\frac{1}{3}(x^3+(1-x)^3 -1)$. (Again, one uses char $(k) \neq
3$.)
\end{remark}

One can check that the functional equation \eqref{3.9} is
equivalent to \eqref{3.4}. To see this, one needs the following
property of the elements $\langle a\rangle$.
\begin{lem}\label{lem3.6} $\langle a\rangle = -a\star\langle
a^{-1}\rangle$.
\end{lem}
\begin{proof}
We remark again that $TB_2(k)
\xrightarrow{\rho\oplus \partial} k\oplus\mathfrak b$ is an
isomorphism, so it suffices to check the relations on $\eps(a)$
and on $\rho(a) = a(1-a)t^2dt$. These become respectively
\ga{3.10}{a\otimes a + (1-a)\otimes (1-a) = a^{-1}\otimes -1 +
(1-a^{-1})
\otimes (1-a) \in k^\times \otimes k \\
-a^3(a^{-1}(1-a^{-1})) = a(1-a). }
The second relation is trivial. For the first one, one writes
\begin{gather}
a\otimes a + (1-a)\otimes (1-a)= a\otimes a + (-a)\otimes (1-a) +
(1-a^{-1})\otimes (1-a)\\
=a^{-1}\otimes (-a+a-1) + (-1)\otimes (1-a).\notag
\end{gather}
Since $k$ is 2-divisible, one has $(-1)\otimes b=0$.
\end{proof}

\section{A conjectural Lie algebra of cycles}

The purpose of this section is to sketch a conjectural algebraic
cycle based theory of additive polylogarithms. The basic
reference is \cite{BK}, where a candidate
for the Tannakian Lie algebra of the category of mixed Tate
motives over a field $k$ is constructed. The basic tool is a
differential graded algebra (DGA) $\sN$ with a supplementary
grading (Adams grading)
\ga{4.1}{\sN^\bullet = \oplus_{j\ge 0} \sN(j)^\bullet \\
\sN(j)^i \subset \text{Codim. $j$ algebraic cycles on }(\P^1-\{1\})^{2j-i}
\notag
}
where $\sN(j)^i$ consists of cycles which meet the faces (defined by setting
coordinates $=0,\infty$) properly and which are alternating with respect to
the action of the symmetric group on the factors and with respect to
inverting the coordinates. The product structure is the external product
$(\P^1-\{1\})^{2j_1-i_1} \times (\P^1-\{1\})^{2j_2-i_2} =
(\P^1-\{1\})^{2j_1-i_1 + 2j_2-i_2}$ followed by alternating projection, and
the boundary map is an alternating sum of restrictions to faces. For full
details, cf. op. cit.

We consider an enlarged DGA
\ga{4.2}{\widetilde\sN^\bullet = \oplus_{j\ge 0} \widetilde\sN(j)^\bullet \\
\widetilde\sN(j)^i := \sN(j)^i \oplus T\sN(j)^i \notag \\
T\sN(j)^i \subset \text{Codim. $j$ algebraic cycles on } \A^1 \times
(\P^1-\{1\})^{2j-i-1}. \notag
}
The same sort of alternation and good position requirements are imposed for
the factors $\P^1-\{1\}$. In addition, we
impose a ``modulus'' condition at
the point $0 \in \A^1$. The following definition is tentative, and is
motivated by example \ref{modex} below.

\begin{defn}\label{defn:mod}Let $D$ be the effective divisor $\A^1\times
\Big((\P^1)^n - \G_m^n\Big)$ on $\A^1\times
\Big((\P^1)^n$, where $\G_m = \P^1-\{0,\infty\}$.
Let $Z\subset\A^1\times (\P^1)^n$ be an effective algebraic cycle. We
assume no component of $Z$ lies on $D$. Let $m\ge 1$ be an
integer. Let $F_i : y_i=1$. We assume no component of $Z$ lies in an
$F_i$. (Components lying in an $F_i$ can be ignored when computing motivic
cohomology).  Write \linebreak $F_i\cdot Z = \sum r_{W,i}W$, and define
$r_W = \max_i(r_{W,i})$. We say that
$Z$ weakly satisfies the modulus $m$ ($Z\equiv 0 \mod m\{0\}\times
(\P^1)^n$) if
the intersection $Z\cdot (0)\times (\P^1)^n = \sum m_V\cdot V$ is defined,
and for each $V$ with $m_V \neq 0$ we have
\eq{4.2a}{m\cdot m_V \le \begin{cases} r_V & V \not\subset D \\
r_V-\epsilon_V & \text{else}\end{cases}
}
Here $\epsilon_V$ is the multiplicity with which $V$ occurs in $Z\cdot D$.

We say $Z$ satisfies modulus $m$ if $Z^0:= Z|_{\A^1\times (\P^1-\{1\})^n}$
is in
good position with respect to all face maps, and if $Z$ and the closures of
all faces of $Z^0$ weakly satisfy modulus $m$.
\end{defn}

If $Z$ satisfies modulus $m$ and $X$ is any subvariety of $(\P^1)^r$ which
is not contained in a face, then $Z\times X$ satisfies modulus $m$ on
$\A^1\times (\P^1)^{n+r}$.

\begin{ex}\label{modex} Milnor $K$-theory of a field can be interpreted in
terms of
$0$-cycles \cite{NS}, \cite{T}. More generally, a Milnor symbol
$\{f_1,\ldots, f_p\}$ over a ring $R$ corresponds to the cycle on
$\Spec(R)\times (\P^1)^p$ which is just the graph
$$\{(x,f_1(x),\ldots f_p(x)) | x\in \Spec(R)\}
$$
One would like cycles with modulus
to relate to
relative $K$-theory. Assume $R$ is semilocal, and let $J\subset R$ be an
ideal. Then we have already used \eqref{2.9} that $K_2(R,J)$ has a
presentation with generators given by pointy-bracket symbols $\langle
a,b\rangle$ with $a\in R$ and $b\in J$ or vice-versa. The pointy-bracket
symbol $\langle a,b\rangle$ corresponds to the Milnor symbol $\{1-ab,b\}$
when the latter is defined. Suppose $R$ is the local ring on $\A^1_k$ at
the origin, with $k$ a field, and take $J= (s^m)$, where $s$ is the
standard parameter. For $a\in J$ and $b\in (s^p)$ for some $p\ge 0$,
 we see that our
definition of cycle with modulus is designed so the cycle
$\{(x,1-a(x)b(x),b(x))\}$ has modulus at least $m$.
\end{ex}

The modulus condition is compatible with pullback to the faces $t_i=0,
\infty$.
\begin{defn}\label{defn4.3} $T\sN(j)^i,\ -\infty\le i\le 2j-1,$ is the
$\Q$-vectorspace of
codimension $j$ algebraic cycles on $\A^1\times(\P^1-\{1\})^{2j-i-1}$ which
are in good position for the face
maps $t_i=0,\infty$ and have modulus $2\{0\}\times (\P^1)^{2j-i-1}$. Here, in
order to calculate the modulus, we close up the cycle to a cycle on
$\A^1\times (\P^1)^{2j-i-1}$. \end{defn}

Note that a cycle $Z$ of modulus $m\ge 1$ doesn't meet $\{0\}\times
(\P^1)^{2j-i-1}$ on  $\A^1\times(\P^1-\{1\})^{2j-i-1}$.

We have a split-exact sequence of $DGA$'s,
\eq{4.4}{0 \to T\sN^\bullet \to \widetilde \sN^\bullet
\stackrel{\longleftarrow}{\to} \sN^\bullet \to 0
}
with multiplication defined so $T\sN^\bullet$ is a square-zero
ideal. Denote the cohomology groups by
\eq{4.5}{\widetilde H^i_M(k,j):= H^i(\widetilde \sN^\bullet(j));\quad
TH^i_M(k,j):= H^i(T\sN^\bullet(j)).
}
As an  example, we will see in section \ref{sec:chow} that the Chow
groups of
$0$-cycles in this context compute the K\"ahler differential forms:
\eq{4.6}{TH^j_M(k,j) \cong \Omega^{j-1}_k; j\ge 0.
}
(Here $\Omega^0_k = k$.)

One may apply the bar construction to the DGA $\widetilde \sN^\bullet$ as
in \cite{BK}. Taking $H^0$ yields an augmented Hopf algebra (defining
$TH^0$ as the augmentation ideal)
\eq{4.7}{0 \to TH^0(B(\widetilde \sN^\bullet)) \to H^0(B(\widetilde
\sN^\bullet)) \stackrel{\longleftarrow}{\to} H^0(\sN^\bullet) \to 0.
}
The hope would be that the corepresentations of the co-Lie algebra of
indecomposables (here $H^{0,+}
:= \ker(H^0 \to \Q)$ denotes the elements of bar degree $>0$,
cf. op. cit. \S 2)
\eq{4.8}{\widetilde \sM := H^0(B(\widetilde
\sN^\bullet))^+/(H^0(B(\widetilde \sN^\bullet))^+)^2 = \sM \oplus T\sM
}
correspond to contravariant motives over $k[t]/(t^2)$. In particular, the
work of Cathelineau \cite{C} suggests a possible additive polylogarithm Lie
algebra. In the remainder of this section, we will speculate a bit on how
this might work.

For a general DGA $A^\bullet$ which is not bounded above, the total grading
on the double complex $B(A^\bullet)$ has infinitely many summands
(cf. \cite{BK}, (2.15)). For example the diagonal line corresponding to
$H^0(B(A^\bullet))$ has terms ($A^+ := \ker(A^\bullet \to \Q$))
\eq{4.9}{A^1,\ (A^+\otimes A^+)^2,\ (A^+\otimes A^+\otimes A^+)^3, \ldots
}
When, however, $A^\bullet$ has a graded structure
\eq{4.10}{A^i = \oplus_{j\ge 0} A^i(j);\quad dA(j) \subset A(j);\quad A^+ =
\oplus_{j>0} A^+(j),
}
for each fixed $j$, only finitely many tensors can occur. For example
$H^0(B(\widetilde \sN^\bullet(1))) = H^1(\widetilde \sN^\bullet(1)) =
k\oplus k^\times$, and $H^0(B(\widetilde \sN^\bullet(2)))$ is the
cohomology along the indicated degree $0$ diagonal in the diagram
\eq{4.11}{\begin{array}{cccc}\widetilde \sN^1(1)\otimes \widetilde \sN^1(1)
& \stackrel{\delta}{\to} & \widetilde\sN^2(2) \\
\uparrow \partial & \makebox[1cm][r]{$\ddots \scriptstyle{\ deg.\ 0}$} &
\uparrow\partial \\
(\widetilde \sN^1(1)\otimes \widetilde \sN^0(1)) \oplus (\widetilde
\sN^0(1)\otimes \widetilde \sN^1(1)) & \stackrel{\delta}{\to} & \widetilde
\sN^1(2) \\
&&\uparrow \\
&& \widetilde\sN^0(2).
\end{array}
}

In the absence of more information about the DGA $\widetilde
\sN^\bullet$, it is difficult to be precise about the indecomposable space
$\widetilde \sM$. As an approximation, we have
\begin{prop} Let $db(\widetilde\sN) \subset \widetilde\sN^1$ be the
subspace
of elements $x$ with decomposable boundary, i.e. such that there exists
$y\in (\widetilde\sN\otimes \widetilde\sN)^2$ with $\delta(y) =
\partial(x)\in \widetilde\sN^2$.  Define
\eq{4.12}{q\widetilde\sN:= db(\widetilde\sN)/(\partial \widetilde\sN^0 +
\delta( \widetilde\sN^+\otimes \widetilde\sN^+)^1)
}
Then there exists a natural map, compatible with the grading by codimension
of cycles (Adams grading)
\eq{4.13}{\phi: \widetilde\sM \to q\widetilde\sN.
}
\end{prop}
\begin{proof}Straightforward.
\end{proof}

As above, we can decompose
\eq{4.14}{q\widetilde\sN = q\sN \oplus Tq\sN,
}
where $q\sN(p)$ is a subquotient of the space of codimension $p$ cycles on
$(\P^1-\{1\})^{2p-1}$, and $Tq\sN$ is a subquotient of the cycles on $\A^1
\times (\P^1-\{1\})^{2p-2}$

\begin{ex}The polylogarithm cycle $\{a\}_p$ for $a\in \C-\{0,1\}$ is
defined to be the image under the alternating projection of
$(-1)^{p(p-1)/2}$ times the locus in $(\P^1-\{1\})^{2p-1}$
parametrized in nonhomogeneous coordinates by
\eq{4.15}{(x_1,\dotsc,x_{p-1},1-x_1,1-x_2/x_1,\dotsc, 1-x_{p-1}/x_{p-2},
1-a/x_{p-1})
}
(We take $\{ a\}_1 = 1-a \in \P^1-\{1\}$.)
To build a class in $H^0(B(\sN^\bullet))^+$ and hence in $\sM$ one uses
that $\partial\{a\}_n = \{a\}_{n-1}\cdot\{1-a\}_1$.
\end{ex}

The following should be compared with \cite{EVG}, where a similar formula
is proposed. The key new point here is that algebraic cycles make it
possible to envision this formula in the context of Lie algebras.
\begin{conj}There exist elements $\langle a\rangle_n \in T\sM(n)$
\eqref{4.8}
represented by cycles $Z_n(a)$ of codimension $n$ on $\A^1\times
(\P^1-\{1\})^{2n-2}$ with $\langle a\rangle_1 = a \in \A^1-\{0\}$. These
cycles should satisfy the boundary condition
\eq{4.16}{\partial \langle a\rangle_n = \langle a\rangle_{n-1}\cdot
\{1-a\}_1 + \langle 1-a\rangle_1\cdot \{a\}_{n-1} \in \Big(\bigwedge^2
\widetilde \sM\Big)(n).
}
\end{conj}

For example, for $n=2$,
\ml{4.17}{\partial \langle a\rangle_2 = a\otimes a + (1-a)\otimes (1-a) \\
 \in
k\otimes k^\times \cong T\sM(1)\otimes \sM(1) \subset \Big(\bigwedge^2
\widetilde \sM\Big)(2)
}
gives Cathelineau's relation \cite{C}.

\begin{prop}Assume given elements $\langle a\rangle_n$ satisfying
\eqref{4.16}. Let $\widetilde \sP = \bigoplus_{n=1}^\infty \Q\langle
a\rangle_n \oplus \Q\{a\}_n$ be the constant graded sheaf over
$\A^1-\{0,1\}$. Then $\widetilde\sP, \partial:\widetilde\sP \to
\bigwedge^2\widetilde\sP$ is a sheaf of co-lie algebras.
\end{prop}
\begin{proof}It suffices to show that $\partial^2=0$. Using the derivation
property of the boundary,
\ml{4.18}{\partial\partial\langle a\rangle_n = (\partial\langle
a\rangle_{n-1})\cdot \{1-a\}_1 - \langle
1-a\rangle_1\cdot\partial\{a\}_{n-1} = \\
\Big(\langle a\rangle_{n-2}\cdot
\{1-a\}_1 + \langle 1-a\rangle_1\cdot \{a\}_{n-2}\Big)\cdot  \{1-a\}_1-
\langle 1-a\rangle_1\cdot\{a\}_{n-2}\cdot\{1-a\}_1 = \\
0 \in \bigwedge^3\widetilde \sM.
}
\end{proof}

We can make the definition (independent of any conjecture)
\begin{defn}The additive polylogarithm sheaf of
Lie algebras over
$\A^1-\{0,1\}$ is the graded sheaf of Lie algebras with graded dual the
sheaf $\widetilde\sP, \partial:\widetilde\sP \to \bigwedge^2\widetilde\sP$
satisfying \eqref{4.16} above.
\end{defn}

Of course, $\langle a\rangle_2$ should be closely related to the element
$\langle a\rangle \in K_2(R,(t^2))$ \eqref{3.1}. The cycle
\eq{4.19}{ \{(t,1-\frac{t^2a(1-a)}{t-1},  \frac{a(1-a)}{t-1})\ |\ t\in
\A^1\} \subset \A^1\times (\P^1-\{1\})^2
}
associated to the pointy bracket symbol in \eqref{3.1} satisfies the
modulus $2$ condition but is not in good position with respect to the
faces. (It contains $(1,\infty,\infty)$.) It is possible to give symbols
equivalent to this one whose corresponding cycle is in good position, but
we do not have a canonical candidate for such a cycle, or a candidate whose
construction would generalize in some obvious way to give all the $\langle
a\rangle_n$.

\section{The Artin-Schreier dilogarithm}
The purpose of this section is to present a definition of what one might
call an Artin-Schreier dilogarithm in characteristic $p$. To begin with, however, we take $X$ to
be a complex-analytic manifold and sketch certain analogies between the multiplicative and
additive theory. We write $\sO$ (resp. $\sO^\times$, $\Omega^1$) for the sheaf of analytic
functions (resp. invertible analytic functions, analytic $1$-forms). The reader is urged to
compare with
\cite{Dsymb}.
\eq{5.1}{\begin{array}{rcl} \text{MULTIPLICATIVE} && \text{ADDITIVE} \\
\sO^\times \stackrel{\mathbb L}{\otimes}
\sO^\times & \longleftrightarrow  & \sO \otimes \sO^\times  \\
 \sO^\times \stackrel{\mathbb L}{\otimes} \sO^\times \to \Big( \sO^\times(1)
\to \Omega^1 \Big)[1] & \longleftrightarrow  & \sO \otimes
 \sO^\times \to \Big( \sO(1)
\to \Omega^1 \Big)[1] \\
K_2 & \longleftrightarrow  & \Omega^1 \\
\text{Steinberg rel'n }= &\longleftrightarrow  &\text{Cathelineau rel'n }= \\
a\otimes (1-a) &  &  a\otimes a + (1-a) \otimes (1-a) \\
\text{exponential of dilogarithm } = &
\longleftrightarrow  &
\text{Shannon entropy function } = \\
\exp\Big(\int_0^a \log(1-t)dt/t\Big) &   & \int^a
\log(\frac{t}{1-t})dt  = \\
& & a\log a + (1-a)\log(1-a).
\end{array}
}

In the multiplicative (resp. additive) theory, one applies $\sO^\times \stackrel{\mathbb
L}{\otimes}_\Z\bullet$ (resp. $\sO\otimes_\Z\bullet$) to the exponential sequence (here $\Z(1) :=
\Z\cdot 2\pi i$)
\eq{5.2}{0 \to \Z(1) \to \sO \to \sO^\times \to 0.
}
The regulator maps \eqref{5.1}, line 2, come from liftings of these tensor
products to
\eq{5.3}{\begin{CD}
\sO^\times(1) @>>> \sO^\times\stackrel{\mathbb L}{\otimes} \sO  @>>>
\sO^\times \otimes \sO^\times\\
@|  @VVV @|\\
\sO^\times(1) @>>> \sO^\times \otimes \sO  @>>> \sO^\times \otimes \sO^\times\\
@VVV @VVV @.\\
\Omega^1 @= \Omega^1.
\end{CD}
}
\eq{5.3b}{\begin{CD}
\sO(1) @>>> \sO \otimes \sO @>>> \sO\otimes \sO^\times\\
@VVV @VVV\\
\Omega^1 @= \Omega^1
\end{CD}
}
In the multiplicative theory, the regulator map can be viewed as associating to two invertible
analytic functions $f,\ g$ on $X$ a line bundle with connection $\sL(f,g)$ on $X$, \cite{Dsymb}.
The exponential of the dilogarithm
\eq{5.4}{\exp(\frac{1}{2\pi i} \int_0^f \log(1-t)\frac{dt}{t})
}
determines a flat section trivializing $\sL(1-g,g)$. Let $\{U_i\}$ be an analytic cover of $X$,
and let $\log_if$ be an analytic branch of the logarithm on $U_i$. Then $\sL(f,g)$ is represented
by the Cech cocycle
\eq{5.5}{\Big(g^{\frac{1}{2\pi i} (\log_i f -\log_j f)}, \frac{1}{2\pi
i}\log_i f\frac{dg}{g}\Big)
}
The trivialization comes from the $0$-cochain
\eq{5.6}{i\mapsto \exp(\frac{1}{2\pi i} \int_0^f \log_i(1-t)\frac{dt}{t}).
}

The additive theory associates to $a\otimes f \in \sO\otimes \sO^\times$
the class in
$\H^1(X, \sO(1) \to \Omega^1)$ represented by the cocycle for $\sO(1) \to
\Omega^1$
\eq{5.7}{\Big(a\otimes (\log_if -\log_j f), \log_i f\cdot da\Big).
}
This can be thought of as defining a connection on the affine bundle
$\sA(a,f)$ associated to the
coboundary of $a\otimes f$ in $H^1(X,\sO(1))$. The affine bundle itself is
canonically trivialized
because in the diagram
\eq{5.8}{\begin{CD}0 @>>> \sO(1) @>>> \sO\otimes \sO @>>> \sO\otimes
\sO^\times @>>> 0 \\
@. @VVV @VVV \\
@. \Omega^1 @= \Omega^1
\end{CD}
}
the top sequence is split (by multiplication $\sO\otimes \sO \to \sO(1)$). The splitting is not
compatible with the vertical arrows, so it does not trivialize the connection. More concretely,
$a\otimes f \in \sO\otimes \sO^\times$ gives the $1$-cocycle $(a\otimes (\log_i f -\log_j f),\
\log_if\cdot da) \in \H^1(X,\sO(1) \to \Omega^1)$. Subtracting the coboundary of the $0$-cochain
$\frac{1}{2\pi i} a\log_i f \otimes 2\pi i$ leaves the cocycle $(0,a\frac{df}{f})$. We have proved:

\begin{prop}\label{prop5.1} The map $\partial : H^0(X,\sO\otimes \sO^\times) \to \H^1(X,\sO(1) \to
\Omega^1)$ factors
$$\begin{CD}H^0(X,\sO\otimes \sO^\times) @>a\otimes f \mapsto adf/f>> H^0(X,\Omega^1) @>>>
\H^1(X,\sO(1) \to \Omega^1).
\end{CD}
$$
In particular, for $a\in \sO$ such that $a$ and $1-a$ are both units, the Cathelineau elements
$\eps(a)=a\otimes a + (1-a)\otimes (1-a)$ \eqref{3.1} lift to
$$a\otimes \log a +(1-a)\otimes \log(1-a) -\frac{1}{2\pi i} \int^a\log(\frac{t}{1-t})dt \otimes
2\pi i \in H^0(X,\sO\otimes \sO).
$$
\end{prop}

\begin{remark}The element $a\otimes a \in H^0(X,\sO\otimes \sO^\times)$
maps to $da=0 \in \H^1(X,\sO(1) \to \Omega^1)$,
but the above construction does not give a canonical trivializing $0$-
cocycle.
\end{remark}

We now suppose $X$ is a smooth variety in characteristic $p>0$, and we
consider an Artin-Schreier analog of the above construction. In place of
the exponential sequence \eqref{5.2} we use the Artin-Schreier sequence of
\'etale sheaves

\eq{5.9}{\begin{CD}0 @>>> \Z/p @>>>  \G_a  @>1-F>>  \G_a @>>> 0.
\end{CD}
}
Here $F$ is the Frobenius map. We replace the twist by
 $\sO^\times_{{\rm an}}$  over $\Z$ by the twist over $\Z/p$ by
$\G_m/\G_m^p$ to build a diagram (compare \eqref{5.3}. Here $Z^1 \subset
\Omega^1$ is the subsheaf of closed forms.)
\eq{5.10}{\begin{CD}0 @>>>  \G_m/\G_m^p @>>>  \G_a
\otimes \G_m/\G_m^p   @>(1-F)\otimes 1>>  \G_a\otimes \G_m/\G_m^p
 @>>> 0 \\
@. @V d\log VV  @VV f\otimes g \mapsto f^pdg/g V \\
@. Z^1 @= Z^1
\end{CD}
}

The group  $\H^1(\G_m \to Z^1)$ is the group of isomorphism classes
of  line bundles with integrable
connections as usual, and $H^1($the subcomplex $\G_m^p \to 0)$
is the subgroup of connections corresponding to a Frobenius descent.
We get an exact sequence
\ml{5.11}{0 \to \{\text{line bundle + integrable connection}\} /\{
\text{lb + Frobenius descent}\} \\
 \to \H^1(X, \G_m/\G_m^p \to Z^1) \to {}_pH^2(X,\G_m) .
}
\begin{prop}Let $\iota, C:Z^1 \to \Omega^1$ be the natural inclusion and
the Cartier operator, respectively. One has a quasi-isomorphism
$(\G_m/\G_m^p \to Z^1) \stackrel{\iota - C}{\sim} \Omega^1[-1]$. Then the
diagram
$$\begin{CD}H^0(X,\G_a\otimes \G_m/\G_m^p) @> \partial\ \eqref{5.10}>>
\H^1(X, \G_m/\G_m^p \to Z^1) \\
@V a\otimes b \mapsto adb/b VV  @V \cong V \iota - C V \\
H^0(X,\Omega^1) @= H^0(X,\Omega^1)
\end{CD}
$$
is commutative.
\end{prop}
\begin{proof} Straightforward from the commutative diagram (with $\mathfrak
b$ defined to make the columns exact and $\phi(a\otimes b) = a\cdot db/b$).
\eq{5.12}{\begin{CD}  @. @. 0 @. 0 \\
 @. @. @VVV @VVV \\
@. @. \mathfrak b @= \mathfrak b \\
@. @. @VVV @VVV \\
0 @>>> \G_m/\G_m^p @>>> \G_a\otimes \G_m/\G_m^p @>F\otimes
1-1>>
\G_a\otimes \G_m/\G_m^p @>>> 0 \\
@. @| @VV\phi\circ(F\otimes 1) V @VV\phi V \\
0 @>>> \G_m/\G_m^p @>>> Z^1 @>1-C>> \Omega^1 @>>> 0 \\
@. @. @VVV @VVV \\
@. @. 0 @. 0.
\end{CD}
}
\end{proof}

Next we want to see what plays the role of the exponential of the
dilogarithm or the Shannon entropy function in this Artin-Schreier
context. Let $X=\Spec(\F_{p^2}[x])$. Begin with Cathelineau's element

\eq{5.13}{\epsilon(x) := x\otimes x + (1-x) \otimes (1-x) \in \mathfrak b.
}
Choose an Artin-Schreier roots $y^p-y=x$ and $\beta^p-\beta=1$. To simplify
we view $\beta \in \F_{p^2}$ as fixed, and we write $\F = \F_{p^2}$. A
local lifting of $\epsilon(x)$ on
the \'etale cover $\Spec\F(y) \to
\Spec\F(x)$ is given by
\eq{5.14}{\rho(y) := y\otimes x + (\beta - y)\otimes (1-x) \in
\Gamma(\Spec\F(y),\G_a\otimes
\G_m/\G_m^p).
}
{From} diagram \eqref{5.12} there should exist a canonical global lifting,
i.e. a
lifting defined over $\Spec\F(x)$. This lifting, call it $\theta(x)$ has
the form $\theta(x) = \rho(y)\cdot \delta(y)^{-1}$ for some $\delta(y) \in
\Gamma(\Spec \F(y),\G_m/\G_m^p)$. We want to calculate $\delta(y)$.

To do this calculation, note
\ml{5.16}{\phi\circ(F\otimes 1)(\rho(y)) = y^pdx/x - (\beta-y)^pdx/(1-x) = \\
\frac{(-y^p(1-y^p+y)+(\beta+1-y^p)(y^p-y))dy}{(y^p-y)(1-y^p+y)} = \\
\frac{(\beta(y^p-y) - y)dy}{(y^p-y)(1-y^p+y)} =:\eta(y) \in Z^1.
}
Viewed as a meromorphic form on $\P^1_y$, $\eta$ has simple poles at the
points $a$ and $\beta -a$ for $a \in \F_p$. The residue of a form $P/Qdy$
at a point $a$ where $Q$ has a simple zero is given by $P(a)/Q'(a)$.
Using this, the residue of $\eta$ at $a\in \F_p$ is $a$. The residue at
$\beta-a$ is $\frac{\beta\cdot 1 -(\beta-a)}{1} = a$. Necessarily,
therefore, since $\eta$ is regular at $y=\infty$ we must have
\eq{5.17}{\eta = d\log\Big(\prod_{a=1}^{p-1}\frac{(\beta-(y+a))^a}{(y+a)^a}
\Big) }
We conclude
\eq{5.18}{\delta(y) = \prod_{a=1}^{p-1}\frac{(\beta-(y+a))^a}{(y+a)^a}.
}
Everything is invariant under the automorphism $y\mapsto
\beta-y$. Indeed, the equation can be rewritten (of course$\mod
\F(y)^{\times p}$)
\eq{5.19a}{\delta(y) =
\prod_{a=1}^{(p-1)/2}\Big[\frac{(\beta-(y+a))(y-a)}{(y+a)(\beta-(y-a))}\Big]^a
}
Note $\delta(y)$ depends on $y$, not just on $x$. Indeed the product
\eqref{6} can be taken for $0\le a\le p-1$, i.e. for $a\in \F_p$. One
gets then

\eq{5.19}{\frac{\delta(y+1)}{\delta(y)} \equiv
\prod_{a=0}^{p-1}\frac{y+a}{\beta-(y+a)} = \frac{y^p-y}{1-y^p-y} =
\frac{x}{1-x}\mod
\F(y)^{\times p}. }
The fact that $\rho(y)\delta(y)^{-1}$ is defined over $\F(x)$ says that the
Cech boundaries of $\rho(y)$ and $\delta(y)$ coincide. Since the latter is,
by definition, the coboundary in $\G_m/\G_m^p$ of $\eps(x) = x\otimes x +
(1-x)\otimes (1-x)$, it follows that $\delta(y)$ is a $0$-cochain for the
Galois cohomology
$$\H^*(\F(y)/\F(x),\G_m/\G_m^p \to Z^1)$$
which trivializes the coboundary of $\rho(\eps(x))$.

Finally, in this section, we discuss a flat realization of the
Artin-Schreier dilogarithm. To see the point, consider the $\ell$-adic
realization of the usual dilogarithm mixed Tate motivic sheaf over
$\A^1-\{0,1\}$. Reducing mod $\ell$ yields a sheaf with fibre an
$\F_\ell$-vector space of dimension $3$. The sheaf has a filtration with
successive quotients having fibres $\Z/\ell \Z, \mu_\ell, \mu_\ell^{\otimes
2}$. The geometric fundamental group acts on the fibre via a Heisenberg
type group. We visualize this action as follows:
\eq{5.20}{\begin{pmatrix} 1 & \mu_\ell & \mu_\ell^{\otimes 2} \\
0 & 1 & \mu_\ell \\
0 & 0 & 1 \end{pmatrix}\begin{pmatrix}\mu_\ell^{\otimes 2} \\
\mu_\ell \\
\Z/\ell\Z \end{pmatrix}
}
Here the notation means that  for $g\in \pi_1^{geo}$ the corresponding matrix
$$\begin{pmatrix} 1 & a_{12}(g) & a_{13}(g) \\
0 & 1 & a_{23}(g) \\
0 & 0 & 1 \end{pmatrix}
$$
has $a_{ij} \in \text{Hom}(\mu_\ell^{\otimes
i-1},\mu_\ell^{\otimes j-1}) = \mu_\ell^{j-i}$.

The essential ingredients here are first the Heisenberg group $\sH_\ell$,
second
the $\sH_\ell$-torsor over $\A^1-\{0,1\}$ corresponding to the kernel of the
representation, and third the (standard) representation of $\sH_\ell$ on
$\Z/\ell\Z \oplus \mu_\ell \oplus \mu_\ell^{\otimes 2}$.

We define an Artin-Schreier Heisenberg group as the non-commutative flat
groupscheme $\sH_{AS}$ which we could suggestively write
\eq{5.21}{\sH_{AS} := \begin{pmatrix}1 & \Z/p\Z & \mu_p \\
0 & 1& \mu_p \\ 0 & 0 & 1\end{pmatrix}.
}
More precisely, $\sH_{AS}$ is
a central extension
\eq{5.26}{0 \to \mu_p \to \sH_{AS} \to \mu_p \times \Z/p\Z \to 0.
}
Let
\ga{5.27}{b : (\mu_p \times \Z/p\Z ) \times (\mu_p \times \Z/p\Z ) \to
\mu_p \\
b((\zeta_1,a_1),(\zeta_2,a_2)) = \zeta_1^{-a_2}\zeta_2^{a_1}. \notag
}
Define $\sH_{AS}=\mu_p \times (\mu_p \times \Z/p\Z)$ as a scheme, with group
structure given by
\eq{5.28}{(\zeta_1,\theta_1,a_1)\cdot (\zeta_2,\theta_2,a_2) :=
(\zeta_1\zeta_2\theta_2^{a_1},\theta_1\theta_2,a_1+a_2).
}
The commutator pairing on $\sH_{AS}$ is given by
\ml{5.29}{\Big[(\zeta_1,\theta_1,a_1), (\zeta_2,\theta_2,a_2)\Big] =
\Big(b((\theta_1,a_1),(\theta_2,a_2)),1,0\Big) = \\
b((\theta_1,a_1),(\theta_2,a_2)) \in \mu_p .}

We fix a solution $\beta^p-\beta = 1$. We define a flat
$\sH_{AS}$-torsor $T = T_\beta$ over $\A^1_{\F_{p^2}}$ as follows. A local (for
the flat topology) section $t$ is determined by \newline
\noindent 1. A $p$-th root of $\frac{x}{1-x}$: $w^p\equiv \frac{x}{1-x}
\mod
\F_{p^2}(x)^{\times p}$.
\newline
\noindent 2. A $y$ satisfying $y^p-y=x$. \newline
\noindent 3. A $p$-th root $z$ of $\delta(y): z^p \equiv \delta(y) \mod
\F_{p^2}(x)^{\times p}$  (where $\delta(y)$ is as in \eqref{5.18}.)

The action of $\sH_{AS}$ is given by
\eq{5.30}{(\zeta,\theta,a)\star(z,w,y) = (\zeta zw^a, \theta
w, y+a).
}
Note $(\zeta zw^a)^p =
\delta(y)(\frac{x}{1-x})^a =
\delta(y+a)$ by
\eqref{5.19}, so the triple on the right lies in $T$. This is an
action because
\ml{5.32}{(\zeta',\theta',a')\star\Big((\zeta,\theta,a)\star(z,w,y)\Big) =
\\ (\zeta',\theta',a')\star\Big((\zeta zw^a, \theta
w, y+a) \Big) = \\
\Big(\zeta'\zeta
zw^a(\theta w)^{a'}, \theta'\theta
w,y+a+a'  \Big) =
(\zeta'\zeta\theta^{a'},\theta'\theta,a+a')\star(z,w,y) = \\
\Big((\zeta',\theta',a')\star(\zeta,\theta,a)\Big)\star(z,w,y).
}

Define $\V := \mu_p \times \mu_p \times \Z/p\Z$. There is an evident action
of $\sH_{AS}$ on $\V$, viewed as column vectors. We suggest that the
contraction $T\stackrel{\sH_{AS}}{\times} \V$ should be thought of as
analogous to the mod $\ell$ \'etale
sheaf on $\A^1-\{0,1\}$ with fibre $\Z/\ell\Z \oplus \mu_\ell \oplus
\mu_\ell^{\otimes 2}$ associated to the $\ell$-adic dilogarithm.

\section{The additive cubical (higher) Chow groups}\label{sec:chow}

In this section, we show that the modulus condition we introduced
in definition \eqref{defn:mod} yields additive Chow groups which
we can compute in weights $(n,n)$. We assume throughout that $k$ is a field
and $\frac{1}{6} \in k$.

One sets
\begin{gather}
A=(\A^1,2\{0\})\\
B=(\P^1\setminus\{1\}, \{0, \infty\}).\notag
\end{gather}
The coordinates will be $x$ on $A$ and $(y_1,\ldots, y_n)$ on $B$.
One considers
\begin{gather}
X_n =A\times B^n.
\end{gather}
The boundary maps $X_{n-1} \inj X_n$ defined by $y_i = 0, \infty$  are
denoted by $\partial_i^j, i=1,\ldots, n, j=0, \infty$.
One denotes by $Y_n\subset X_n$ the union of the faces
$\partial_i^j(X_{n-1})$.
One defines
\begin{gather}
\sZ_0(X_n)=\oplus \Z \xi, \ \xi \in X_n \setminus Y_n, \\
\xi \ \text{closed point}. \notag
\end{gather}
For any $1$-cycle  $C$ in $X_n$, one denotes by
 $\nu: \bar{C}\to \P^1\times (\P^1)^n$ the
normalisation of its compactification.
One defines
\begin{gather}
\sZ_{1}(X_n)=\oplus \Z C, \ C\subset X_n \ \text{with}\\
\partial_i^j(C) \in \sZ_0(X_{n-1}) \
\text{and (cf. definition \ref{defn:mod})} \notag\\
2\nu^{-1}(\{0\}\times (\P^1)^n)+ \nu^{-1}(Y_n)
\subset {\rm max}_{i=1}^n
\nu^{-1}(\P^1\times (\P^1)^{i-1} \times \{1\} \times (\P^1)^{n-i})\notag
\notag
\end{gather}
One defines
\begin{gather}
\partial:=\sum_{i=1}^n (-1)^{i} (\partial_i^0-\partial_i^\infty):
\sZ_{1}(X_{n})\to \sZ_0(X_{n-1}) \ \text{for \ all} \ i, j.
\end{gather}
Further one  defines the differential form
\begin{gather}
\psi_n =
\frac{1}{x} \frac{dy_1}{y_1} \wedge \ldots \wedge \frac{dy_n}{y_n}
\in
\Gamma(\P^1\times (\P^1)^n,
\Omega^1_{(\P^1\times (\P^1)^n)/\Z} (\log Y_n)(\{x=0\})).
\end{gather}
We motivate the choice of this differential form  as follows.
One considers
\begin{gather}
V_n(t)=(\P^1\setminus \{0, t\}, \infty)\times (\P^1\setminus\{0, \infty\},1)^n.
\end{gather}
Its cohomology
\begin{gather}
H^{n+1}(V_n(t))= H^1\otimes (H^1)^n=F^1\otimes (F^1)^n
\end{gather}
is Hodge-Tate for $t\neq 0$. The  generator is given by
\begin{gather}
\omega_{n+1}(t)= (\frac{dx}{x} -\frac{d(x-t)}{(x-t)})\wedge \frac{dy_1}{y_1}
\wedge \ldots \wedge \frac{dy_n}{y_n}.
\end{gather}
Thus
\begin{gather}
\frac{\omega_{n+1}(t)}{t}|_{t=0}=d(\psi_n).
\end{gather}
\begin{defn}\label{defn:addchow}
We define the additive cubical (higher) Chow groups
\eq{}{ TH_{M}^n(k,n)= \sZ_0(X_{n-1})/\partial \sZ_{1}(X_n)\notag}
\end{defn}
One has the following reciprocity law
\begin{prop} \label{recip}
The map $\sZ_0(X_{n-1})\to \Omega^{n-1}_k$ which associates to a closed
point $\xi\in X_{n-1}\setminus Y_{n-1}$ the value
 ${\rm Trace}(\kappa(\xi)/k)
(\psi_{n-1}(\xi))$ factors through $$TH_{M}^n(k,n):=\sZ_0(X_{n-1})
/\partial\sZ_1(X_{n}).$$
\end{prop}
\begin{proof}
Let $C$ be in $\sZ_1(X_{n})$. Let $\Sigma\subset \bar{C}$ be the locus
of poles of $\nu^*\psi_{n}$. One has the functoriality map
\begin{gather}
\nu^*: \Omega^n_{\P^1\times (\P^1)^{n-1}}(\log Y_{n-1})
(\{x=0\})\to \Omega^{n-1}_{\bar{C}/\Z}(*\Sigma).
\end{gather}
Thus reciprocity says
\begin{gather}
\sum_{\sigma \in \Sigma} {\rm res}_\sigma \nu^*(\psi_{n})=0.
\end{gather}
Recall that here res means the following. One has
a surjection
\begin{gather}
\Omega^n_{\bar{C}/\Z}(*\Sigma) \to
\Omega^{n-1}_{k/\Z}\otimes \omega_{\bar{C}/k}(*\Sigma)
\end{gather}
which yields
\begin{gather}
\Gamma(\bar{C}, \Omega^n_{\bar{C}/\Z}(*\Sigma))
\to \Gamma(\bar{C}, \Omega^{n-1}_{k/\Z}\otimes \omega_{\bar{C}/k}(*\Sigma))=
\Omega^{n-1}_{k/\Z} \otimes \Gamma(\bar{C}, \omega_{\bar{C}/k}(*\Sigma)).
\end{gather}
 By definition, res on $\Omega^{n-1}_{k/\Z}\otimes
\omega_{\bar{C}/k}(*\Sigma)$ is $1\otimes {\rm res}$.
This explains the reciprocity.

Now we analyze $\Sigma\subset \sigma^{-1}(Y_n \cup \{x=0\})$.
Let $t$ be a local parameter on $\bar{C}$ in a point $\sigma$ of
$\nu^{-1}(\{x=0\})$.

We write
$x=t^m\cdot u$, where $u\in \sO^\times_{\bar{C}, \sigma}, m\ge 0$.
If $m\ge 1$, the assumption we have on $\sZ_1$ says that there is at least
 one $i$ such that
$\{t=0\}$ lies in $\nu^{-1} (\{y_i=1\})$. Let us order $i=1,\ldots, n$
such that
$\{t=0\}$ lies in $\nu^{-1} (\{y_i=1\})$ for $i=1,2,\ldots, r$. Thus we
write
\begin{gather}
y_i-1=t^{m_i}\cdot u_i, m_1 \ge m_2\ge \ldots
 \ge 1, u_i \in \sO^\times, i=1,\\
 y_i=t^{p_i}u_i, p_i\ge 0, u_i\in \sO^\times, i=r+1,\ldots, n.\notag
\end{gather}
The assumption we have says
\begin{gather}
2m\le m_1.
\end{gather}
One has around the point $\sigma$
\begin{gather}
\nu^{-1}(\psi_n)|_{\sigma}= \frac{u^{-1}}{t^m}\cdot
\frac{d(t^{m_1}\cdot u_1)}{1+t^{m_1}\cdot u_1} \wedge \ldots \wedge
\frac{d(t^{m_r}\cdot u_r)}{1+t^{m_r}\cdot u_r}\wedge_{i=r+1}^n
 \frac{d(t^{p_i}\cdot u_i)}{t^{p_i}\cdot u_i}.
\end{gather}
We analyze the poles of the right hand side. The numerator of this expression
is divisible by $t^{(m_1+\ldots + m_r)-1}$.
Thus the condition for
$\nu^{-1}\psi_n$ to be smooth in $\sigma$ is
\begin{gather}
m+1\le (m_1+\ldots +m_r).
\end{gather}
This is always fulfilled for $2m\le m_1$. We have
\begin{gather}
\nu^{-1}\psi_n \ \text{smooth \ in } \ \nu^{-1} (\{x=0\}).
\end{gather}
On the other hand, one obviously has
\begin{gather}
{\rm res}_{y_i=0}\psi_n=-{\rm res}_{y_i=\infty}\psi_n=(-1)^i \psi_{n-1}.
\end{gather}
Thus one concludes
\begin{gather}
\sum_{\sigma \in \Sigma} {\rm res}_\sigma \nu^{-1}\psi_n= \sum_{i=1}^n (-1)^i
\psi_{n-1} (\partial_i^0-\partial_i^\infty)(C)=0.
\end{gather}
\end{proof}
We want to see that the reciprocity map in proposition \ref{recip} is an
isomorphism.
Define
\begin{gather}\label{6.23}
k\otimes_{\Z} \wedge_{i=1}^{n-1} k^\times \to TH_{M}^n(k,n)\\
a\otimes (b_1\wedge \ldots \wedge b_{n-1}) \mapsto
 (\frac{1}{a}, b_1,\ldots, b_{n-1}) \ \text{for} \ a\neq 0 \notag\\
\mapsto 0 \ \text{for} \ a=0.\notag
\end{gather}
\begin{prop} Assume $\frac{1}{6} \in k$. Then \eqref{6.23} factors through
$$\Omega^{n-1}_k \to TH_{M}^n(k,n).$$
\end{prop}
\begin{proof}
One has the following relations, where $\equiv$ means equivalence modulo
$\partial\sZ_1(X_{n+1})$:
\begin{gather}
\label{6.26} (\frac{1}{x+x'}, y_1,\ldots,y_n)\equiv
(\frac{1}{x},y_1,\ldots, y_n) +
(\frac{1}{x'}, y_1,\ldots, y_n)\\
\label{6.27}(x, y_1z_1,y_2,\ldots, y_n)\equiv (x,y_1,y_2,\ldots, y_n) +
(x, z_1, y_2,\ldots, y_n) \\
(x, -1, y_2,\ldots, y_n)\equiv 0 \in TH_{M}^n(k,n).
\end{gather}
Note, the last is obviously a consequence of the first two:
\begin{gather}
(2x,-1,  y_2, \ldots, y_n)=2(x, -1, y_2,\ldots, y_n)=(x, 1, y_2,\ldots,
y_n)= 0,
\end{gather}
so we need only consider \eqref{6.26} and \eqref{6.27}. Assume first
$xx'(x+x') \neq 0$, and define
\begin{gather}
C=(t, y_1=\frac{(1-xt)(1-x't)}{1-(x+x')t}, y_2,\ldots, y_n)\in
\sZ_1(X_n).
\end{gather}
Indeed, the expansion of $y_1$ in
$t=0$ reads $ 1+t^2c_2 + ($higher order terms$)$,
so our modulus condition is
fulfilled. Also we have taken $y_i \in k^\times$ so $C$
 meets the faces properly.  Then one has
\begin{gather}
\partial(C)= (\frac{1}{x}, y_2,\ldots, y_n) + (\frac{1}{x'},
y_2,\ldots, y_n) - (\frac{1}{x+x'}, y_2, \ldots, y_n).
\end{gather}
Similarly, if $x+x'=0$, then one sets
\begin{gather}
C=(t, y_1=(1-\frac{t^2}{x^2}), y_i)\in \sZ_1(X_n).
\end{gather}
One has
\begin{gather}
\partial(C)= (x, y_i) + (-x, y_i),
\end{gather}
proving \eqref{6.26}. Note the proposition for $n=1$ is a consequence of
this identity.

To show multiplicativity in the
$y$ variables, one uses Totaro's  curve \cite{T}. There is a $\sC\in
\sZ_1(B^{n+1})$ with $\partial(\sC)=(y_1z_1, y_2, \ldots,
y_n)-(y_1, y_2, \ldots, y_n) -(z_1, y_2,\ldots, y_n)$. One sets
$C=(x, \sC) \in \sZ_1(X_{n+1})$. Here $x$ is fixed and nonzero, so the
modulus condition is automatic, and one has $\partial(C)=(x,
\partial(\sC))$. This proves \eqref{6.27}.

It remains to verify the Cathelineau relation (cf. \cite{BE}, \cite{C})
\begin{gather}
(\frac{1}{a}, a,b_2, \ldots, b_{n})+
(\frac{1}{1-a}, (1-a), b_2,\ldots, b_{n}) \equiv 0.
\end{gather}
In fact, the $b_2,\ldots b_n \in k^\times$ play no role, so we will drop
them. One considers the 1-cycle which is given by its parametrization
\begin{gather}
Z(a)=-Z_1(a)+Z_2\\
 Z_1(a)=(t, 1+\frac{t}{2}, 1-\frac{a^2t^2}{4})\notag\\
Z_2= ( \frac{t}{4}, 1+ \frac{t}{6}, 1-\frac{t^2}{4}), \notag
\end{gather}
We see immediately that $Z\in \sZ_1(X_{2})$.
One has
\begin{gather}
\partial(Z_1(a))\\
=(-2, 1-a^2)-(\frac{2}{a},
1+\frac{1}{a})-(-\frac{2}{a}, 1-\frac{1}{a})=\notag \\
(-2, 1-a)  + (-2, 1+a) + (\frac{2}{a}, \frac{a-1}{a+1})= \notag \\
(-2, 1-a) + (\frac{2}{a}, a-1) + (-2, 1+a)
-(\frac{2}{a}, a+1)=\notag\\
(-2, a-1) + (\frac{2}{a}, a-1) + (-2, 1+a) +
(-\frac{2}{a}, a+1)=\notag\\
(\frac{2}{a-1},a-1) -(\frac{2}{a+1}, a+1).\notag
\end{gather}
Setting $a=1-2b$, one obtains
\begin{gather}
\partial(Z_1(a))= (-\frac{1}{b}, -2b)-(\frac{1}{1-b}, 2(1-b))=\\
-(\frac{1}{b}, b)-(\frac{1}{1-b}, 1-b) -(\frac{1}{b},
2) -(\frac{1}{1-b}, 2)=\notag \\
-(\frac{1}{b}, b)-(\frac{1}{1-b}, 1-b) -(1, 2).\notag
\end{gather}
One has
\begin{gather}
\partial(Z_2)= (-\frac{3}{2}, -8) -(\frac{1}{2}, \frac{4}{3})
-(-\frac{1}{2}, \frac{2}{3})=\\
3(-\frac{3}{2}, 2) -(\frac{1}{2}, 2)=\notag\\
(-\frac{3}{2}, 2) + (-\frac{3}{4}, 2)-(\frac{1}{2}, 2)=\notag\\
(-\frac{1}{2}, 2) -(\frac{1}{2}, 2)= -(1, 2).\notag
\end{gather}
In conclusion
\begin{gather}
\partial(Z(1-2a))=(\frac{1}{a}, a, b_i) + (\frac{1}{1-a}, 1-a,
b_i).
\end{gather}
\end{proof}

We now have well defined maps $\phi_n, \psi_n$
\begin{gather}
\Omega^{n-1}_k \xrightarrow{\phi_n} TH_{M}^n(k, n) \xrightarrow{\psi_n}
\Omega^{n-1}_k
\end{gather}
which split $TH_M^n(k, n)$. The image of the differential forms
consists of all 0-cycles which are equivalent to 0-cycles $\sum
m_i p_i$ with $p_i\in X_n(k)$.
\begin{thm} . Assume $\frac{1}{6} \in k$. The above maps identify
$TH_{M}^n(k, n)$ with
$\Omega^{n-1}_k$.
\end{thm}
\begin{proof}
It suffices to show that the class of a give closed point $p \in
(\A^1-\{0\}) \times
(\P^1-\{0, 1, \infty\})^{n-1}$ lies in the image of $\phi_n$. Write
$\kappa = \kappa(p)$ for the residue field at $p$.
One first applies a Bertini type argument as in \cite{BE}, Proposition 4.5,
to reduce to the case where $\kappa/k$ is separable. Then we follow the
argument in loc.cit. The degree
$[\kappa:k] < \infty$, so standard cycle constructions yield a norm map
$N:TH_M^n(\kappa,n) \to TH_M^n(k,n)$. We claim the diagram
\eq{6.39}{\begin{CD}\Omega^{n-1}_\kappa @>\phi_{n,\kappa}>>
TH_M^n(\kappa,n) \\
@V \tr VV @VV N V \\
\Omega^{n-1}_k @> \phi_{n,k} >> TH_M^n(k,n)
\end{CD}
}
is commutative, where $\tr$ is the trace on differential forms. Indeed,
$\Omega^n_{\kappa}=\kappa\otimes \Omega^n_k$, so it
suffices to check on forms $ad\log(b_1)\wedge \ldots \wedge d\log(b_{n-1})$
with $a\in \kappa$ and $b_i \in k$. But in this situation, we have
projection formulas, both for $0$-cycles and for differential forms, and it
is straightforward to check (ignore the $b_i$ and reduce to $n=1$)
\ga{6.41}{\tr(ad\log(b_1)\wedge \ldots \wedge d\log(b_{n-1})) =
\tr_{\kappa/k} (a)d\log(b_1)\wedge \ldots \wedge d\log(b_{n-1}) \\
N(\frac{1}{a},b_1,\ldots,b_{n-1}) =
\begin{cases}(\frac{1}{\tr_{\kappa/k}a},b_1,\ldots,b_{n-1}) & \tr(a) \neq 0
\\
0 & \tr(a) = 0 \end{cases} \notag
}
Write $[p]_\kappa$ (resp. $[p]_k$) for the class of $p \in
TH^n_M(\kappa,n)$ (resp. $TH^n_M(k,n)$.) One has
$[p]_k=N([p]_\kappa)$. Since $p$ is $\kappa$-rational, $[p]_\kappa \in
\text{Image}(\phi_\kappa)$. Commutativity of \eqref{6.39} implies $[p]_k
\in \text{Image}(\phi_k)$. It follows that $\phi_k$ is surjective, proving
the theorem. \end{proof}

\bibliographystyle{plain}
\renewcommand\refname{References}

\end{document}